\documentclass[12pt]{amsart}
\usepackage{amsmath}
\usepackage{amsxtra}
\usepackage{amscd}
\usepackage{amsthm}
\usepackage{amsfonts}
\usepackage{amssymb}
\usepackage{eucal}

\input prepictex
\input pictex
\input postpictex

\newtheorem{thm}{Theorem}[section]
\newtheorem{lem}{Lemma}[section]
\newtheorem{cor}{Corollary}[section]
\newtheorem{prop}{Proposition}[section]

\theoremstyle{definition}

\theoremstyle{remark}
\newtheorem{rem}{Remark}[section]

\numberwithin{equation}{section}

\begin{document}

\newcommand{\thmref}[1]{Theorem~\ref{#1}}
\newcommand{\secref}[1]{Section~\ref{#1}}
\newcommand{\lemref}[1]{Lemma~\ref{#1}}
\newcommand{\propref}[1]{Proposition~\ref{#1}}
\newcommand{\corref}[1]{Corollary~\ref{#1}}
\newcommand{\remref}[1]{Remark~\ref{#1}}
\newcommand{\eqnref}[1]{(\ref{#1})}
\newcommand{\exref}[1]{Example~\ref{#1}}

\newcommand{\nc}{\newcommand}
 \nc{\on}{\operatorname}
 \nc{\Z}{{\mathbb Z}}
 \nc{\C}{{\mathbb C}}
 \nc{\oo}{{\mf O}}
 \nc{\R}{{\mathbb R}}
 \nc{\N}{{\mathbb N}}
 \nc{\bib}{\bibitem}
 \nc{\pa}{\partial}
 \nc{\F}{{\mf F}}
 \nc{\rarr}{\rightarrow}
 \nc{\larr}{\longrightarrow}
 \nc{\al}{\alpha}
 \nc{\ri}{\rangle}
 \nc{\lef}{\langle}
 \nc{\W}{{\mc W}}
 \nc{\gam}{\ol{\gamma}}
 \nc{\Q}{\ol{Q}}
 \nc{\q}{\widetilde{Q}}
 \nc{\la}{\lambda}
 \nc{\ep}{\epsilon}
 \nc{\g}{\mf g}
 \nc{\B}{\mf B}
 \nc{\h}{\mf h}
 \nc{\Hy}{ \widetilde{H} }
 \nc{\n}{\mf n}
 \nc{\A}{{\mf a}}
 \nc{\G}{{\mf g}}
 \nc{\HH}{{\mf h}}
 \nc{\Li}{{\mc L}}
 \nc{\La}{\Lambda}
 \nc{\is}{{\mathbf i}}
 \nc{\V}{\mf V}
 \nc{\bi}{\bibitem}
 \nc{\NS}{\mf N}
 \nc{\dt}{\mathord{\hbox{${\frac{d}{d t}}$}}}
 \nc{\E}{\mc E}
 \nc{\ba}{\tilde{\pa}}
 \def\smapdown#1{\big\downarrow\rlap{$\vcenter{\hbox{$\scriptstyle#1$}}$}}
 \nc{\hf}{\frac{1}{2}}
 \nc{\mc}{\mathcal}
 \nc{\mf}{\mathfrak}
 \nc{\ol}{\fracline}
 \nc{\el}{\ell}
 \nc{\etabf}{{\bf \eta}}
 \nc{\x}{{\bf x}}
 \nc{\xibf}{{\bf \xi}}
 \nc{\y}{{\bf y}}
 \nc{\NP}{\Pi}
 \nc{\gltwo}{{\rm gl}_{\infty|\infty}}
 \nc{\hz}{\hf+\Z}
 \nc{\hsd}{\widehat{\mc S\mc D}}
 \nc{\parth}{\partial_{\theta}}
 \nc{\sd}{\mc S \mc D}

 \nc{\D}{\mathcal D}
 \nc{\gl}{{\rm gl}_\infty}
 \nc{\hD}{\widehat{\mathcal D}}
 \nc{\hgl}{\widehat{\rm gl}_{\infty}}
 \nc{\vac}{|0 \rangle}
 \nc{\winfty}{\mathcal W_{1+\infty}}
 \nc{\hgltwo}{\widehat{\rm gl}_{\infty}}

\title{The Bloch-Okounkov correlation functions at higher levels}
\author[Shun-Jen Cheng]{Shun-Jen Cheng}
\address{Department of Mathematics, National Taiwan University, Taipei,
Taiwan 106} \email{chengsj@math.ntu.edu.tw}

\author[Weiqiang Wang]{Weiqiang Wang}
\address{Department of Mathematics, University of Virginia, Charlottesville, VA 22904}
\email{ww9c@virginia.edu}

\begin{abstract}
We establish an explicit formula for the $n$-point correlation
functions in the sense of Bloch-Okounkov for the irreducible
representations of $\hgl$ and $\mathcal W_{1+\infty}$ of arbitrary
positive integral levels.
\end{abstract}

\keywords{Correlation functions, Fock space, $\hgl$, $\winfty$
algebra.}

\subjclass[2000]{Primary: 17B65; Secondary: 17B68.}

\maketitle

\section{Introduction}\label{intro}

The representation theories of the infinite-dimensional  Lie
algebras $\hgl$ and $\winfty$ have been well developed. Among the
most interesting representations are the integrable highest weight
modules of a positive integral level (cf.~\cite{AFMO, DJKM, FKRW,
KR} and the references therein). In \cite{BO}, Bloch and Okounkov
formulated certain $n$-point correlation functions on the
fermionic Fock space and established a beautiful closed formula
for them in terms of theta functions. In light of free field
realizations of $\hgl$ and $\winfty$, their work amounts to the
study of certain trace functions on the irreducible modules of
$\hgl$ or $\winfty$ of level $1$, which can be regarded as a
character formula of these modules involving all elements in the
infinite-dimensional Cartan subalgebras of $\hgl$ or $\winfty$.

The aim of this article is to formulate and establish an explicit
formula for the $n$-point correlation function for an arbitrary
integrable highest weight module of $\hgl$ or $\winfty$ of any
positive integral level $\ell\in\N$. Our main result is
\thmref{mainthm}, which gives an elegant simple formula for these
correlation functions essentially as the product of $\el$ copies
of the ``normalized'' Bloch-Okounkov $n$-point function of level
$1$ and the $q$-dimension formula of the corresponding integrable
module of level $\el$. The main tool used here is a duality
between finite-dimensional representations of ${\rm GL}(\ell)$ and
the integrable modules of $\hgl$ or $\winfty$ of level $\ell$ (cf.
\cite{F, FKRW, W1}). This duality is an infinite-dimensional
generalization of the Howe duality \cite{Howe}. Using the
Bloch-Okounkov formula \thmref{mainthm} is then derived by
exploiting several combinatorial consequences of this duality. We
also present another formula for these $n$-point correlation
functions of higher levels involving inverse Kostka numbers. As an
immediate consequence some combinatorial identities are obtained
by comparing these two formulas.

Besides the obvious representation-theoretic motivations, the
interest in these $n$-point functions is stimulated by the fact
that Bloch-Okounkov's correlation functions of level $1$ also have
connections to combinatorics of partitions (cf. \cite{BO, Ok}) and
geometry (such as Gromov-Witten theory of an elliptic curve by
Okounkov-Pandharipande, equivariant intersections on Hilbert
schemes by Li-Qin-Wang). It would be interesting to see whether
these $n$-point functions at higher levels afford similar
geometric applications. We also refer to \cite{Mil, W2} for a
twisted analog of \cite{BO} and interesting connections to vertex
operators.

The paper is organized as follows. In
Section~\ref{sec:formulation} we introduce the necessary notation
and formulate the problem of $n$-point correlation functions of
higher levels. We also recall the Bloch-Okounkov formula and
conclude the section by stating our main theorem.
Section~\ref{sec:main} is divided into three subsections and its
main purpose is to provide a proof of the main theorem.  In
\secref{sec:kostka} we derive another formula for these $n$-point
correlation functions of higher levels involving inverse Kostka
numbers and from it some combinatorial identities.

\section{Formulation and statement of the main theorem} \label{sec:formulation}

Recall that any finite-dimensional irreducible rational
representation of the general linear group ${\rm GL}(\ell)$ is
obtained by taking the tensor product of an integral power of the
determinant representation with an irreducible polynomial
representation. Thus, the irreducible rational representations of
${\rm GL}(\ell)$ are parameterized by {\em generalized partitions}
of length $\ell$, where by a generalized partition $\la$ of length
$\ell$ we mean a non-increasing sequence of integers
$(\la_1,\la_2,\cdots,\la_\ell)$. For a generalized partition $\la$
of length $\ell$, we will denote by $|\la|=\sum_{i=1}^\ell\la_i$
the {\em size} of $\la$. Given such a $\la$, we denote the
corresponding irreducible ${\rm GL}(\ell)$-module by $V_\ell^\la$.

Let $\gl$ denote the Lie algebra spanned by matrices of the form
$(a_{ij})_{i,j\in\Z}$ with $a_{ij}=0$, for $|i-j|>>0$, that is,
${\rm gl}_\infty$ is the Lie algebra of infinite matrices with
finitely many non-zero diagonals.  The Lie algebra
$\gl=\sum_{j\in\Z}(\gl)_j$ is $\Z$-graded by setting
$\text{deg}(E_{ij})=j-i$, where $E_{ij}$ denotes the matrix with
$1$ in the $i$th row and $j$th column and zero elsewhere.
Furthermore it has a non-trivial $2$-cocycle $\alpha$ defined by
$\alpha(A,B)={\rm Tr}([J,A]B)$, which gives rise to a central
extension $\hgl$ \cite{DJKM}. Here $J=\sum_{i\le 0}E_{ii}$. We
shall denote by $K$ the central element of this central extension.
The Lie algebra $\hgl$ inherits from $\gl$ the $\Z$-gradation
$\hgl=\sum_{j\in\Z}(\hgl)_j$ with $\text{deg} (K) =0$, and its
Cartan subalgebra is $(\hgl)_0=(\gl)_0+\C K$. Thus any element
$\La$ in the restricted dual $(\hgl)^*_0$ gives rise to a highest
weight irreducible $\hgl$-module $L(\hgl,\La)$ of highest weight
$\La$.

The fundamental weights of $\hgl$ are given as follows. For $i>0$
we set $\La_i(K)=1$ and
\begin{equation*}
\La_i(E_{jj})=\left\{\begin{array}{l l} 1,& {\rm if}\ 0<j\le
i,\\0, & {\rm otherwise}.
\end{array}\right.
\end{equation*}
Now if $i< 0$ we set $\Lambda_i(K)=1$ and
\begin{equation*}
\La_i(E_{jj})=\left\{\begin{array}{l l} -1,& {\rm if}\ i<j\le
0,\\0, & {\rm otherwise}.
\end{array}\right.
\end{equation*}
Finally we define $\La_0$ by $\La_0(K)=1$ and $\La_0(E_{jj})=0$,
for all $j\in\Z$.

Given a positive integer $\ell$ and a generalized partition
$\la=(\la_1,\cdots,\la_\ell)$ of length $\ell$, we denote by
$\La(\la)$ the $\hgl$-highest weight
$\La_{\la_1}+\cdots+\La_{\la_\ell}$. The highest weight modules
$L(\hgl,\La(\la))$ with highest weight vector $v_{\La(\la)}$ of
level $\el$ are the so-called integrable modules, which are
arguably the most interesting $\hgl$-modules.

We define
\begin{eqnarray}\label{operatorT}
 T(t) = \sum_{k \in \Z} t^{k-\hf} E_{kk} +
 \frac{1}{t^{\hf}-t^{-\hf}}\; K,
\end{eqnarray}
where $t$ is an indeterminate.

Extending the definition of Bloch and Okounkov \cite{BO} of level
$1$, we define the level $\ell$ {\em $n$-point correlation
function} on $L(\hgl,\La(\la))$ by
\begin{equation*}
F_\la^\ell (q;t_1,\cdots,t_n) :={\rm
Tr}_{L(\hgl,\La(\la))}(q^{H}T(t_1)T(t_2)\cdots T(t_n)).
\end{equation*}
Here $H$ is the energy operator characterized by
\begin{eqnarray}
H \cdot v_{\La(\la)} &=& (\la^2/2) \cdot v_{\La(\la)}, \label{weight} \\
 {[}H, E_{ij}] &=& (i-j) E_{ij}, \nonumber
\end{eqnarray}
where by definition $$\la^2 :=\la_1^2 +\la_2^2 +\cdots +
\la_\el^2.$$ The \eqnref{weight} is introduced for convenience
later on (see \eqnref{energy} below).

A remarkable closed formula for the $n$-point function at level
$1$ has been established in \cite{BO}.  Before stating it we need
some additional notation.

Let $q$ denote a formal parameter. Let $\Theta(t)$ denote the
following theta function
\begin{equation*}
\Theta(t) =\Theta_{11}(t;q) :=
\sum_{n\in\Z}(-1)^nq^{\frac{(n+\hf)^2}{2}}t^{n+\hf}.
\end{equation*}
It is agreed that $1/ (-k)!=0$ for $k>0$ below. Let
\begin{equation*}
\Theta^{(k)}(t) :=(t\frac{d}{dt})^k\Theta(t),
\end{equation*}
and let $\varphi(q)$ denote the Euler product
\begin{equation*}
\varphi(q)=\prod_{j=1}^\infty(1-q^j).
\end{equation*}

\begin{thm}\label{BO} \text{\rm (Bloch-Okounkov \cite{BO})}
The correlation function $F^1_{(0)} (q; t_1,\cdots,t_n)$ of level
$1$ equals
\begin{eqnarray}  \label{eq:bo}
\frac{1}{\varphi(q)}\cdot \sum_{\sigma\in S_n} \frac{{\rm det}
\Big{(}\frac{\Theta^{(j-i+1)}(t_{\sigma(1)}\cdots
t_{\sigma(n-j)})}{(j-i+1)!} \Big{)}_{i,j=1}^n}
{\Theta(t_{\sigma(1)}) \Theta(t_{\sigma(1)}t_{\sigma(2)})\cdots
\Theta(t_{\sigma(1)}t_{\sigma(2)}\cdots t_{\sigma(n)})}.
\end{eqnarray}
\end{thm}
The expression \eqnref{eq:bo} will be denoted by
$F_{BO}(q;t_1,\cdots,t_n)$.  It is implicitly given in \cite{Ok}
that, for $k \in \Z$,
$$F^1_{(k)}(q;t_1,\cdots,t_n)
=q^{\frac{k^2}{2}}(t_1\cdots t_n)^kF_{BO}(q;t_1,\cdots,t_n).$$

Our main result is the following formula for the Bloch-Okounkov
$n$-point correlation functions at level $\el$.  The proof will be
given in \secref{sec:main}.

\begin{thm}\label{mainthm}
The $n$-point correlation function of level $\ell$ associated to
the generalized partition $\la=(\la_1, \cdots, \la_\ell)$ is given
by
\begin{align*}
F^\ell_\la(q; t_1,\cdots,t_n)=
 &\frac{q^{\frac{\la^2}{2}}(t_1t_2\cdots t_n)^{|\la|}{\prod_{1\le i<j\le
\el}\Big{(}1-q^{\la_i-\la_j+j-i}\Big{)}}}{\varphi(q)^\el} \times\\
&\left (\sum_{\sigma\in S_n} \frac{{\rm det}
\Big{(}\frac{\Theta^{(j-i+1)}(t_{\sigma(1)}\cdots
t_{\sigma(n-j)})}{(j-i+1)!} \Big{)}_{i,j=1}^n}
{\Theta(t_{\sigma(1)}) \Theta(t_{\sigma(1)}t_{\sigma(2)})\cdots
\Theta(t_{\sigma(1)}t_{\sigma(2)}\cdots
t_{\sigma(n)})}\right)^\el.
\end{align*}
\end{thm}

\begin{rem}  \label{rem:winfty}
The $\winfty$ algebra $\hD$ is the central extension of the Lie
algebra $\D$ of differential operators on the circle. Write $D = t
\frac{d}{dt}$ and elements in $\D$ can be written as a linear
combination of $t^k f_k(D)$, $k \in \Z$, where the $f_k$'s are
polynomials in one variable.
                %
%
Setting the degrees of $t^{-k}, D$ and $C$ to $k, 0$ and $0$,
respectively defines a $\Z$-grading of $\hD$. A Lie algebra
homomorphism $\phi$ of $\hD$ into $\hgl$, preserving the
$\Z$-gradation, is given by
$$
 \phi \left( t^k f( D) \right) = \sum_{j \in \Z} f (-j )  E_{j-k, j}.
$$
The pullback of the integrable module $L(\hgl, \La)$ via $\phi$
remains an irreducible $\hD$-module, and the most interesting
$\hD$-modules are obtained in this way or its variant (cf.
\cite{FKRW, KR}). In this way, the $n$-point functions introduced
above can be regarded as $n$-point functions for irreducible
$\hD$-modules of level $\el$.
\end{rem}

\section{Proof of the main theorem} \label{sec:main}

\subsection{The $({\rm GL}(\el),\hgl)$-Duality on $\F^\el$}

Let $z$ be a formal indeterminate and, for $1\le i\le \ell$, let
\begin{eqnarray*}
\psi^{+i}(z) &=& \sum_{r\in\hf+\Z}\psi^{+i}_r z^{-r-\hf}, \\
\psi^{-i}(z) &=& \sum_{r\in\hf+\Z}\psi^{-i}_r z^{-r-\hf}
\end{eqnarray*}
denote $\ell$ pairs of free fermionic fields whose non-trivial
(super)commutation relations are given by:
\begin{equation*}
[\psi^{+i}_r,\psi^{-j}_s]=\delta_{ij}\delta_{r,-s} \text{I},\quad
1\le i,j\le \ell.
\end{equation*}
Let $\F^\ell$ denote the Fock space associated to these $\ell$
pairs of fermions with vacuum vector $\vac$ satisfying
\begin{eqnarray*}
\psi^{\pm i}_r \vac =0, \qquad 1 \le i \le \ell, \quad r>0.
\end{eqnarray*}

In terms of the generating function $ E(z,w) =\sum_{i,j \in \Z}
E_{ij} z^{i -1} w^{ -j}$,
we have a free field realization for $\hgl$ of level $\ell$ (i.e.
the central element $K$ acts as $\ell \cdot \text{I}$) on
$\F^\ell$ given by:
 \begin{eqnarray}\label{gl:gen}
   E (z, w) & =& \sum_{i=1}^\ell :\psi^{+i} (z) \psi^{-i} (w):.
 \end{eqnarray}
The normal ordering $:\ :$ here and further is defined by moving
the annihilation operators to the right (up to a sign).

Put
\begin{equation*}
e_{ij}=\sum_{r\in\hf+ \Z}:\psi_{r}^{+i}\psi^{-j}_{-r}:, \quad 1\le
i,j\le \ell.
\end{equation*}
Then the $e_{ij}$'s define a representation of the general linear
Lie algebra ${\rm gl}(\ell)$ on $\F^\ell$, which lifts to an
action of the Lie group ${{\rm GL}}(\ell)$.  In particular as a ${
{\rm GL}}(\ell)$-module $\F^{\ell}$ is completely reducible.

\begin{prop} \label{duality} \cite{F} {\rm (}cf.~\cite{W1}{\rm )}
The pair $({\rm GL}(l),\hgl)$ on $\F^\ell$ forms a dual reductive
pair in the sense of Howe.
Furthermore, as a ${\rm GL}(\ell)\times\hgl$-module, $\F^\ell$ is
multiplicity-free and decomposes into
\begin{equation} \label{eq:dual}
\F^\ell=\bigoplus_{\la} V_\ell^\la \bigotimes L(\hgl,\La(\la)),
\end{equation}
where $\la$ is summed over all generalized partitions of length
$\el$.
\end{prop}

\begin{rem}
As explained in \cite{FKRW}, the duality \eqnref{duality} extends
to a $({\rm GL}(\ell), \winfty)$ duality via the homomorphism in
Remark~\ref{rem:winfty}.
\end{rem}

The energy operator $H$ on $\F^\el$ can be realized as
\begin{equation}\label{energy}
H=\sum_{i=1}^\el\sum_{r\in\hf+\Z} r :\psi^{+i}_{-r}
\psi^{-i}_{r}:,
\end{equation}
which may be regarded as the zero-mode of a Virasoro field.

Similar to the definition of monomial symmetric polynomials and
Schur polynomials, one may define the monomial symmetric (Laurent)
polynomial $m_\la$ and Schur (Laurent) polynomial $s_\la$
associated to a generalized partition $\la$. The $s_\la
(z_1,\ldots, z_\ell)$'s, where the $\la$'s are generalized
partitions of length $\ell$, form a linear basis for $\C[z_1^{\pm
1}, \ldots, z_\ell^{\pm \ell}]^{S_\ell}$.  Similar statement is
true for the $m_\la (z_1,\ldots, z_\ell)$'s.

\begin{cor}\label{combid}
We have the following combinatorial identity.
\begin{align}  \label{eq:comb}
\prod_{i=1}^\el\prod_{r\in\hf+\Z_+}(1+q^rz_i)(1+q^rz_i^{-1})& =\\
\sum_{\la}s_{\la}(z_1,\cdots,z_l)&\frac{q^{\frac{\la^2}{2}}{\prod_{1\le
i<j\le
\el}\Big{(}1-q^{\la_i-\la_j+j-i}\Big{)}}}{{\varphi(q)^\el}}\nonumber,
\end{align}
where $\la=(\la_1,\cdots,\la_\ell)$ is summed over all generalized
partitions of length $\ell$.
\end{cor}

\begin{proof}

We apply the trace of the operator
$q^Hz_1^{e_{11}}z_2^{e_{22}}\cdots z_\el^{e_{\el\el}}$ to both
sides of the duality \eqnref{eq:dual}. Applying it to the
left-hand side of \eqnref{eq:dual} gives rise to
\begin{equation*}
{\rm Tr}_{\F^\el} (q^Hz_1^{e_{11}}z_2^{e_{22}}\cdots
z_\el^{e_{\el\el}}) =
\prod_{i=1}^\el\prod_{r\in\hf+\Z_+}(1+q^rz_i)(1+q^rz_i^{-1}).
\end{equation*}

On the other hand applying the operator
$q^Hz_1^{e_{11}}z_2^{e_{22}}\cdots z_\el^{e_{\el\el}}$ to the
right-hand side of \eqnref{eq:dual} we see that
$z_1^{e_{11}}z_2^{e_{22}}\cdots z_\el^{e_{\el\el}}$ acts on the
first tensor factor, while $q^H$ acts on the second tensor factor.
Now the $q$-dimension formula for the $\hgl$-module
$L(\hgl,\La(\la))$ is well-known to be (cf. e.g. \cite{FKRW})
\begin{equation*}
{\rm Tr}_{L(\hgl,\La(\la))} q^H =
\frac{q^{\frac{\la^2}{2}}\prod_{1\le i<j\le
\el}\Big{(}1-q^{\la_i-\la_j+j-i}\Big{)}}{\varphi(q)^\el},
\end{equation*}
while $s_\la (z_1,\ldots, z_\ell)$ is the character of
$V_\ell^\la$. Thus the trace operator when applied to the
right-hand side of \eqnref{eq:dual} gives rise to the right-hand
side of \eqnref{eq:comb}.
\end{proof}

\subsection{Traces of operators on $\F^\el$}
It follows from \eqnref{gl:gen} that the operator $T(t)$ in
\eqnref{operatorT} acting on $\F^\ell$ is given by the formula
\begin{equation*}
T(t) = \sum_{i=1}^\ell \sum_{r\in\hf+\Z} t^{r} \psi^{+i}_{-r}
\psi^{-i}_{r}.
\end{equation*}
Let $S$ be the {\em shift operator} on $\F^1$, which is uniquely
determined by the relations
\begin{eqnarray*}
 S^{-1}\psi^{\pm}_r S
  &=& \psi^{\pm}_{r\pm 1}, \quad r\in\hf+\Z, \\
 S \, \vac
  &=& \psi_{-\hf}^{+} \vac.
\end{eqnarray*}
For $\el =1$, $e_{11} \in {\rm GL}(1)$ coincides with the standard
charge operator $C$ on $\F^1$. Introduce

$${\bf F}(z,q;t_1,\cdots,t_n)={\rm Tr}_{\F^1} (z^{e_{11}} q^H
T(t_1)\cdots T(t_n)).$$
 We have the following lemma (compare with \cite{Ok}).

\begin{lem}\label{auxlem}
We have
\begin{equation*}
{\bf F}(z,q;t_1,\cdots,t_n)= F_{BO}(q;t_1,\cdots,t_n)
\sum_{k\in\Z}(zt_1t_2\cdots t_n)^k q^{\frac{k^2}{2}}.
\end{equation*}
\end{lem}

\begin{proof}
\propref{duality} for $\ell =1$ is simply the boson-fermion
correspondence: $\F^1=\sum_{k\in\Z} V^{k}_1\otimes
L(\hgl,\La(k))$, where $V^k_1$ is the $1$-dimensional module over
${\rm GL}(1)$ of (highest) weight $k$. The space $L(\hgl,\La(k))$
is exactly the  charge $k$ subspace $\F_{(k)}$ of $\F^1$. The
following commutation relations are standard (cf. e.g. \cite{Ok}):
\begin{eqnarray*}
  S^k \F_{(0)} = \F_{(k)}, &&\quad
  S^{-1}T(t) S =t\, T(t) \\
  S^{-1} C S = C+1, &&\quad
  S^{-1}H S =H+C+\hf.
\end{eqnarray*}
Using these relations, we compute that
\begin{align*}
{\bf F}(z,q; t_1,\cdots,t_n)
=& \sum_{k\in\Z}{\rm Tr}_{\F_{(k)}}
(z^C q^H T(t_1)\cdots T(t_n)) \allowdisplaybreaks\\
=&\sum_{k\in\Z}{\rm Tr}_{\F_{(0)}} (S^{-k}z^C q^HT(t_1)\cdots
T(t_n)S^{k})\allowdisplaybreaks\\
=&\sum_{k\in\Z}{\rm Tr}_{\F_{(0)}} (z^{C+k} q^{H+\frac{k^2}{2}+kC}
T(t_1)\cdots
T(t_n))(t_1\cdots t_n)^k\allowdisplaybreaks\\
=&\sum_{k\in\Z}{\rm Tr}_{\F_{(0)}} ((zt_1\cdots t_n)^k
q^{H+\frac{k^2}{2}}
T(t_1)\cdots T(t_n)) \allowdisplaybreaks\\
=&\sum_{k\in\Z}(zt_1\cdots t_n)^k q^{\frac{k^2}{2}}{\rm
Tr}_{\F_{(0)}}(q^H
T(t_1)\cdots T(t_n))\allowdisplaybreaks\\
=&F_{BO}(q;t_1,\cdots,t_n)\sum_{k\in\Z}(zt_1\cdots
t_n)^kq^{\frac{k^2}{2}}.
\end{align*}
This finishes the proof.
\end{proof}

\begin{lem}  \label{lem:identity}
We have the following combinatorial identity:
\begin{equation}\label{auxmain2}
\prod_{i=1}^\ell{\bf
F}(z_i,q;t_1,\cdots,t_n)=\sum_{\la}s_\la(z_1,\cdots,z_\ell)F^\ell_\la(q;t_1,\cdots,t_n),
\end{equation}
where $\la$ above is summed over all generalized partitions of
length $\el$.
\end{lem}

\begin{proof}
The lemma follows from the computation of the trace
$${\rm Tr}_{\F^\ell} (z_1^{e_{11}}\cdots z_\ell^{e_{\ell\ell}} q^H
T(t_1)\cdots T(t_n))$$
in two different ways. Note that all the operators involved here
commute with one another.

First note that $\F^\ell$ is the tensor product of the $\ell$ Fock
spaces associated to each of the pairs of fermions $\psi^{\pm
i}(z)$, $1 \le i \le \ell$. The operators $q^H$ and $T(t_i)$ act
diagonally on this tensor product, while $z^{e_{ii}}$ acts
nontrivially only on the $i$-th tensor factor. These
considerations and the definition of $\bf F$ give rise to the
left-hand side of \eqnref{auxmain2}.

On the other hand, we may compute the trace by using the Howe
duality \eqnref{eq:dual}. The $z_1^{e_{11}}\cdots
z_\ell^{e_{\ell\ell}}$ acts on the first tensor factor
$V^\la_\ell$ only and $q^H T(t_1)\cdots T(t_n)$ acts only on the
second factor $L(\hgl,\La(\la))$. Noting that ${\rm Tr}
(z_1^{e_{11}}\cdots z_\ell^{e_{\ell\ell}})$ on $V_\ell^\la$ is
exactly $s_\la(z_1, \ldots, z_\ell)$, we obtain the right-hand
side of \eqnref{auxmain2}.
\end{proof}

\subsection{Completion of the proof of \thmref{mainthm}}

It follows from \lemref{auxlem} that
\begin{equation}\label{aux10}
\prod_{i=1}^\ell{\bf F}(z_i,q; t_1,\cdots,t_n)
=F_{BO}(q;t_1,\cdots,t_n)^\ell \prod_{i=1}^\el\sum_{k_i\in\Z}
(z_it_1\cdots t_n)^{k_i}q^{\frac{k_i^2}{2}}.
\end{equation}
Using the celebrated Jacobi triple product identity
\begin{equation*}
\prod_{j=1}^\infty(1-q^j)(1+x_iq^{j-\hf})(1+x_i^{-1}q^{j-\hf}) =
\sum_{k\in\Z}x_i^{k_i}q^{\frac{k_i^2}{2}},
\end{equation*}
we can rewrite \eqnref{aux10}, with $x_i=z_it_1\cdots t_n$, as
\begin{align*}
\prod_{i=1}^\ell{\bf F}&(z_i,q; t_1,\cdots,t_n) =\\
&F_{BO}(q;t_1,\cdots,t_n)^\ell\varphi(q)^\el
\prod_{i=1}^\el\prod_{j=1}^\infty(1+(z_it_1\cdots
t_n)q^{j-\hf})(1+(z_it_1\cdots t_n)^{-1}q^{j-\hf}).
\end{align*}
Using \corref{combid}, now with $z_it_1\cdots t_n$ replacing
$z_i$, we thus obtain
\begin{align}\label{aux21}
&\prod_{i=1}^\ell{\bf F}(z_i,q; t_1,\cdots,t_n)
=\\
&\sum_{\la}s_{\la}(z_1t_1\cdots t_n,\cdots,z_\el t_1\cdots
t_n)F_{BO}(q;t_1,\cdots,t_n)^\ell \cdot q^{
\frac{\la^2}{2}}{\prod_{1\le i<j\le
\el}\Big{(}1-q^{\la_i-\la_j+j-i}\Big{)}}=\nonumber\\
&\sum_{\la}s_{\la}(z_1,\cdots,z_\el )(t_1\cdots
t_n)^{|\la|}F_{BO}(q;t_1,\cdots,t_n)^\ell \cdot q^{
\frac{\la^2}{2}}{\prod_{1\le i<j\le
\el}\Big{(}1-q^{\la_i-\la_j+j-i}\Big{)}}.\nonumber
\end{align}
The last identity follows from the homogeneity of the Schur
polynomials. Noting that the Schur functions associated to
generalized partitions are linearly independent, we can now
complete the proof of \thmref{mainthm} by comparing the
coefficient of $s_\la(z_1,\cdots,z_l)$ in \eqnref{aux21} with the
one in \eqnref{auxmain2}. \qed

\begin{cor}\label{trivialpar}
The $n$-point correlation function of the vacuum module of $\hgl$
{\rm (}i.e. associated to the trivial partition $(0)${\rm )} at
level $\el$ is given by
\begin{equation*}
F^\el_{(0)}(q;t_1,\cdots,t_n)=F_{BO}(q;t_1,\cdots,t_n)^\ell
{\prod_{1\le i<j\le \el}\Big{(}1-q^{j-i}\Big{)}}.
\end{equation*}
\end{cor}

The next corollary follows from \corref{combid} and the Jacobi
triple product identity.

\begin{cor}
We have the following combinatorial identity:
\begin{equation*}
\prod_{i=1}^\el\sum_{k_i\in\Z}z_i^{k_i}q^{\frac{k_i^2}{2}}=\sum_{\la}
s_\la(z_1,\cdots,z_\el)q^{\frac{\la^2}{2}}{{\prod_{1\le i<j\le
\el}\Big{(}1-q^{\la_i-\la_j+j-i}\Big{)}}},
\end{equation*}
where $\la=(\la_1,\cdots,\la_\ell)$ is summed over all generalized
partitions of length $\ell$.
\end{cor}

\begin{rem}
In \cite{W1} several other types of Howe duality, realized on some
Fock spaces, were established. They typically involve pairs
consisting of a finite-dimensional Lie group of classical type and
an infinite-dimensional Lie algebra, which is a classical Lie
subalgebra of $\hgl$. We expect that the duality method used in
this paper can be extended to compute the analogous $n$-point
functions for the irreducible modules of these Lie algebras (cf.
\cite{W2} for results in a very special case). This direction will
be pursued elsewhere.
\end{rem}

\section{Another formula for the $n$-point
functions}\label{sec:kostka}

Recall that the {\em Kostka number} $K_{\la\mu}$ (cf.~\cite{Mac})
associated to two partitions $\la$ and $\mu$ is determined by
\begin{equation}\label{kostka}
s_\la=\sum_{\mu}K_{\la\mu} m_\mu.
\end{equation}
Since the monomial symmetric Laurent polynomials and the Schur
Laurent functions corresponding to generalized partitions are
linear bases for $\C[z_1^{\pm 1}, \ldots, z_\ell^{\pm
\ell}]^{S_\ell}$, we can define the Kostka number $K_{\la\mu}$ by
means of \eqnref{kostka}, for generalized partitions $\la$ and
$\mu$.

Let $r$ be a non-negative integer such that $\la+(r^\ell)$ and
$\mu+(r^\ell)$ are partitions, where $\la+(r^\ell)$ denotes
$(\la_1+r, \cdots, \la_\ell +r)$. One clearly has
\begin{eqnarray}
m_{\la+(r^\ell)} =&(z_1\cdots z_\ell)^r m_{\la},\quad
s_{\la+(r^\ell)} =(z_1 \cdots z_\ell)^r s_{\la},\nonumber\\
&K_{\la+(r^\ell),\mu+(r^\ell)} = K_{\la\mu}.\label{shift}
\end{eqnarray}
In particular $K_{\la\mu}$ is zero unless both $\la$ and $\mu$
have the same size.

\begin{thm}\label{aux100} The $n$-point correlation function associated to the
generalized partition $\la$ of level $\ell$ is given by
\begin{equation}
F^\ell_\la (q; t_1,\cdots,t_n)=(t_1t_2\cdots t_n)^{|\la|}F_{BO}
(q;t_1,\cdots,t_n)^\ell \cdot \sum_\mu
q^{\frac{\mu^2}{2}}K^{(-1)}_{\mu\la},
\end{equation}
where the summation above is over all generalized partitions $\mu$
of length $\ell$, and $(K^{(-1)}_{\mu\la})$ denotes the inverse
Kostka matrix.
\end{thm}

\begin{proof}
It follows from \lemref{auxlem} that
\begin{align} \label{aux0}
\prod_{i=1}^\ell{\bf F}&(z_i,q; t_1,\cdots,t_n)
=\\
&F_{BO}(q;t_1,\cdots,t_n)^\ell \sum_{(k_1,\cdots,k_\ell)
\in\Z^\ell} (t_1\cdots t_n)^{k_1+\cdots+k_\el}z_1^{k_1} \cdots
z_\ell^{k_\ell} q^{\sum_{i=1}^\ell \frac{{k_i}^2}{2}}.\nonumber
\end{align}

We may rewrite
\begin{equation}\label{aux1}
\sum_{(k_1,\cdots,k_\ell)\in\Z^\ell}(t_1\cdots
t_n)^{k_1+\cdots+k_\el}z_1^{k_1} \cdots z_\ell^{k_\ell}
q^{\frac{{k_i}^2}{2}}=\sum_{\mu}(t_1\cdots t_n)^{|\mu|}m_\mu (z_1,
\ldots, z_\ell) q^{\frac{\mu^2}{2}},
\end{equation}
summed over the generalized partition $\mu$ of length $\ell$.

By recalling the definition of the Kostka numbers $K_{\la\mu}$
\eqnref{kostka}, we may rewrite the right-hand side of
\eqnref{aux1} as
\begin{equation} \label{aux2}
\sum_{\mu}(t_1\cdots t_n)^{|\mu|}m_\mu(z_1, \ldots, z_\ell)
q^{\frac{\mu^2}{2}} =\sum_{\la,\mu, |\la|=|\mu|}  (t_1\cdots
t_n)^{|\mu|}q^{\frac{\mu^2}{2}} K^{(-1)}_{\mu\la}
s_\la(z_1,\cdots,z_\ell).
\end{equation}

By combining \eqnref{aux0}, \eqnref{aux1} and \eqnref{aux2}, we
arrive at
\begin{align}\label{auxmain1}
\prod_{i=1}^\ell {\bf F}&(z_i,q; t_1,\cdots,t_n)
=\\
&\sum_{\la}(t_1\cdots t_n)^{|\la|}F_{BO}(q; t_1,\cdots,t_n)^\ell
\Big{(} \sum_{\mu,|\mu|=|\la|} q^{\frac{\mu^2}{2}}
K^{(-1)}_{\mu\la}\Big{)} s_\la (z_1,\cdots,z_\ell).\nonumber
\end{align}

Now the Schur functions $s_\la (z_1,\cdots,z_\ell)$'s are linearly
independent. Hence the theorem follows from comparing
\eqnref{auxmain2} with \eqnref{auxmain1}.
\end{proof}

{From} \thmref{mainthm} and \thmref{aux100} we obtain the
following.

\begin{cor}\label{aux1000}
For a generalized partition $\la=(\la_1,\cdots,\la_\el)$ of length
$\el$ we have
\begin{equation}\label{combid3}
\sum_\mu q^{\frac{\mu^2}{2}}K^{(-1)}_{\mu\la} =
q^{\frac{\la^2}{2}}{\prod_{1\le i<j\le
\el}\Big{(}1-q^{\la_i-\la_j+j-i}\Big{)}},
\end{equation}
where the summation $\mu$ is over generalized partitions of the
same size as $\la$.
\end{cor}

In the case of ordinary partitions \corref{aux1000} may be
restated as follows.

\begin{cor}
Given a partition $\la$ we have the following combinatorial
identity:
\begin{equation}\label{combid2}
\sum_{\mu}K^{(-1)}_{\mu\la}q^{\frac{\mu^2}{2}}+\sum_{r=1}^{\el-1}\sum_{l(\mu)<\el}K^{(-1)}_{\mu,\la+(r^\el)}
q^{\frac{\sum_{i=1}^\ell
(\mu_i-r)^2}{2}}=q^{\frac{\la^2}{2}}{\prod_{1\le i<j\le
\el}\Big{(}1-q^{\la_i-\la_j+j-i}\Big{)}},
\end{equation}
where the $\mu$'s are partitions and $l(\mu)$ denotes the length
of $\mu$. In particular if $\la$ is the trivial partition, we
obtain
$$
\sum_{r=0}^{\el-1} \sum_{l(\mu)<\el}K^{(-1)}_{\mu(r^\el)}
q^{\frac{\sum_{i=1}^\ell (\mu_i-r)^2}{2}}=\prod_{1\le i<j\le
\el}(1-q^{j-i}).
$$
\end{cor}

\begin{proof}
The proof amounts to showing that the left-hand side of
\eqnref{combid3} equals the left-hand side of \eqnref{combid2}.
This follows from the combinatorial formula (2) on page 107 in
\cite{Mac} for $K^{(-1)}_{\mu\la}$ and \eqnref{shift} (note that
we have interchanged the roles of $\la$ and $\mu$ from {\em loc.
cit.}). We remark that an ingredient in the proof is the fact that
$K_{\mu\la}=0$ for $\la_\el\ge \el$ and $l(\mu)<\el$, which is a
consequence of the combinatorial formula for $K^{(-1)}_{\mu\la}$
in {\em loc. cit.}.
\end{proof}

\bigskip

\noindent{\bf Acknowledgements.} The first author is partially
supported by an NSC-grant of the R.O.C. He also wishes to express
his gratitude to the Institute of Mathematical Science and the
Department of Mathematics of the University of Virginia for
support and hospitality. We thank Antun Milas for very interesting
and stimulating correspondence. We were informed that he has
independently studied the problem addressed in this paper along
the lines of \cite{Mil}.


\bigskip
\frenchspacing

\begin{thebibliography}{AB}

\bibitem{AFMO} H. Awata, M. Fukuma, Y. Matsuo, and S. Odake,
{\em Representation theory of the $\winfty$ algebra}, Prog. Theor.
Phys. Suppl. {\bf 118} (1995), 343--373.

\bi{BO} S.~Bloch and A.~Okounkov, {\em The character of the
infinite wedge representation}, Adv.~Math.~{\bf 149} (2000),
1--60.

\bibitem{DJKM} E.~Date, M.~Jimbo, M.~Kashiwara and T.~Miwa,
{\em Operator approach to the Kadomtsev-Petviashvili equation.
Transformation groups for soliton equations III}, J. Phys. Soc.
Japan {\bf 50} (1981), 3806--3812.

\bi{F} I.~Frenkel, {\em Representations of affine Lie algebras,
Hecke modular forms and Kortweg-de Vries type equations},
Lect.~Notes.~Math.~{\bf 933} (1982), 71--110.

\bi{Howe} R.~Howe, {\it Perspectives on invariant theory: Schur
duality, multiplicity-free actions and beyond}, The Schur
Lectures, Israel Math. Conf. Proc. {\bf 8}, Tel Aviv (1992),
1--182.

\bi{FKRW} E.~Frenkel, V.~Kac, A.~Radul, and W.~Wang, {\em
$W_{1+\infty}$ and $W(gl_N)$ with central charge $N$}, Commun.
Math. Phys. {\bf 170} (1995), 337--357.


\bi{KR} V.~Kac and A.~Radul, {\em Quasi-finite highest weight
modules over the Lie algebra of differential operators on the
circle}, Commun. Math. Phys. {\bf 157} (1993), 429--457.

\bi{Mac} I.~G.~Macdonald, {\em Symmetric functions and Hall
polynomials}, Second Edition, Oxford Math.~Monogr., Clarendon
Press, Oxford, 1995.

\bibitem{Mil} A.~Milas,
{\em Formal differential operators, vertex operator algebras and
zeta-values, II}, J.~Pure Appl.~Algebra {\bf 183} (2003),
191--244.

\bibitem{Ok} A.~Okounkov,
{\em Infinite wedge and random partitions}, Selecta~Math.~New
Series {\bf 7} (2001), 1--25.


\bi{W1} W.~Wang, {\em Duality in infinite dimensional Fock
representations}, Commun.~Contem.~Math. {\bf 1} (1999), 155--199.

\bi{W2} W.~Wang, {\em Correlation functions of strict partitions
and twisted Fock spaces}, Transform.~Groups (to appear),
math.QA/0303259.



\end{thebibliography}
\end{document}